\documentstyle[12pt]{article}
\newcommand{\nc}{\newcommand}

\newcommand{\Hom}{\,{\rm Hom}\,}

\newcommand{\beq}{\begin{equation}}
\newcommand{\eeq}{\end{equation}}
\newcommand{\beqst}{\begin{equation*}}
\newcommand{\eeqst}{\end{equation*}}
\newcommand{\barr}{\begin{array}}
\newcommand{\earr}{\end{array}}
\newcommand{\beqar}{\begin{eqnarray}}
\newcommand{\eeqar}{\end{eqnarray}}
\newtheorem{theorem}{Theorem}[section]
\newtheorem{lemma}[theorem]{Lemma}
\newtheorem{prop}[theorem]{Proposition}
\newtheorem{definition}[theorem]{Definition}
\newtheorem{remit}[theorem]{Remark}

\newenvironment{rem}{\begin{remit}\rm}{\end{remit}}

\newcommand{\RR}{{{\bf  R }}}
\newcommand{\CC}{{{\bf  C }}}
\nc{\FF}{ {\Bbb F} } 
\newcommand{\ZZ}{{{\bf   Z }}}

\newcommand{\UU}{{\Bbb U }}

\newcommand{\cald}{{\mbox{$\cal D$}}}
\newcommand{\cale}{{\mbox{$\cal E$}}}

\newcommand{\call}{{\mbox{$\cal L$}}}
\newcommand{\calm}{{\mbox{$\cal M$}}}

\newcommand{\calo}{{\mbox{$\cal O$}}}

\newcommand{\calr}{{\mbox{$\cal R$}}}
\newcommand{\cals}{{\mbox{$\cal S$}}}
\newcommand{\calt}{{\mbox{$\cal T$}}}

\def\g{\gamma}

\def\l{\lambda}

\begin{document}
\nc{\Proof}{ \noindent{\bf Proof:} }
\nc{\ctwo}{{F^c}}
\nc{\qtwo}{{F^q}}
\nc{\orb}[1]{{ \calo_{#1} } } 
\nc{\redpar}{{t}}
\nc{\altorpar}{\Lambda}
\nc{\orpar}{{\Lambda}}
\nc{\cone}{Q}
\nc{\proj}{\Pi}
\nc{\push}[1]{\iota_{#1} }
\nc{\diffop}{{\call}}
\nc{\dr}[1]{D_{#1}}
\nc{\Conv}{{ \rm Conv} }
\nc{\cthree}{ {S^c}}
\nc{\liek}{{ \bf {g}}}
\nc{\liets}{\liet^*}
\nc{\lieks}{\liek^*}
\nc{\lietpl}{{ \liet_+}}
\nc{\lietplo}{{ \liet^o_+}}
\nc{\lietspl}{\liets_+}
\nc{\lietsplo}{(\liets_+)^o }
\nc{\modclass}{{\cals^c} }
\nc{\qthree}{{\volsppar{0}{3}}}
\nc{\modquant}{{ \cals^q}}
\nc{\threeorpar}{{ \orpar_1, \orpar_2, \orpar_3 }}
\nc{\threealtorpar}{{ \altorpar_1, \altorpar_2, \altorpar_3 }}
\nc{\threeredpar}{{ \redpar_1, \redpar_2, \redpar_3 }}
\nc{\group}{{G}}
\nc{\torus}{{ T}}
\nc{\domp}{ {D_+} }
\nc{\Del}{\bigtriangleup} 
\nc{\afun}{a}
\nc{\atil}{\tilde{a} } 
\nc{\charac}[1]{\epsilon(#1)}
\nc{\intlat}{\Lambda^I}
\nc{\winv}{{ \frac{1}{|W|} }}
\nc{\conjclass}{ {\rm Cl} } 
\nc{\modsp}{ { \calm} }
\nc{\modsppar}[2]{\modsp_{#1,#2} } 
\nc{\volsppar}[2]{S_{#1,#2} } 
\nc{\volparone}{\volsppar{g}{1}(\Lambda)}     
\nc{\volpar}{\volsppar{g}{n} (\bt) } 
\nc{\modspparreg}[2]{{\modsp^o_{#1,#2} } }
\nc{\upmodsppar}[2]{\calr_{#1,#2} } 
\nc{\upmodspparreg}[2]{\calr^o_{#1,#2} } 
\nc{\sig}[2]{\Sigma^{#1}_{#2} }
\nc{\tc}[1]{[{#1}] }
\nc{\mommap}[1]{\mu^{#1}}
\nc{\mommapo}{{\mu}}
\nc{\modpar}{\modsppar{g}{b} (\Lambda^{(1)}, \dots, \Lambda^{(b)}) } 
\nc{\modparc}{\modsppar{g}{b} (c, \bfl) } 
\nc{\modparone}{\modsppar{g}{1} (\Lambda) } 
\nc{\modparonec}{\modsppar{g}{1} (c,\Lambda) } 
\nc{\extmodsppar}[2]{\tilde{\calm}_{#1,#2} } 
\nc{\extmodpar}{\extmodsppar{g}{n} } 
\nc{\extmodparone}{\extmodsppar{g}{1} } 
\nc{\inprd}[1]{{ (#1) }}
\nc{\inprt}[1]{{ \langle #1 \rangle }}
\nc{\curve}[1]{{ C_{#1} }}
\nc{\momexto}{{J}}
\nc{\momext}[1]{{\momext_{#1} }}
\nc{\alc}{{ D_+} }
\nc{\alco}{{ D^o_+} }
\nc{\Ad}{ {\rm Ad}}
\nc{\bt}{ {\bf { \Lambda}} }
\nc{\bd}{ {c} }
\nc{\mapsing}{{ \phi }}
\nc{\tbunsppar}[3]{{ V_{#1,#2}^{(#3)}  } }
\nc{\tbunspparcirc}[4]{{ V_{#1,#2,#3}^{(#4)}  } }
\nc{\tbunpar}{\tbunsppar{g}{n}{m} (\bt) }
\nc{\tbunparone}{\tbunsppar{g}{1}{m} (\Lambda) }
\nc{\tbunparnn}{\tbunsppar{g}{n}{1} (\bt) }
\nc{\tbunparcirc}{\tbunspparcirc{g}{n}{\alpha}{m} (\bt) }
\nc{\tbunparcircg}{\tbunspparcirc{g}{n}{\gamma_j}{m} (\bt) }
\nc{\loopp}[1]{{ [S_{#1} ]    }}
\nc{\tilal}{\tilde{\alpha} }
\nc{\lietp}{{ \liet^\perp} }
\nc{\nplus}{{ n_+} }
\nc{\dist}{{ f} }
\nc{\qfour}{{ S_{0,4,f } }}
\nc{\fouraltorpar}{{ \altorpar_1, \altorpar_2, \altorpar_3, \altorpar_4 } }
\nc{\waff}{ { W_{\rm aff} }}
\nc{\torpar}{{ h}}
\nc{\tps}{{ P }}
\nc{\qtw}{{ q}}
\nc{\dimg}{{d}}
\nc{\dimt}{{ l}}
\nc{\pant}{ {P}}
\nc{\triup}{ { \bigtriangleup} } 
\nc{\free}{ { \Bbb F } } 
\nc{\univ}{ {\Bbb U} } 
\nc{\flagm}{ {\rm Fl} }
\nc{\hattau}{ { \hat{\tau} } }
\nc{\hu}[1]{ {\hat{u}_{#1} } }
\nc{\modonetwo}[1]{ \calm_{1,1}^{{SU(2)}} (#1) }
\nc{\tortwo}{ { T_2 } }
\nc{\tbunonetwo}[1]{ V_{1,1}^{{SU(2) }} (#1) }
\nc{\subv}{ { Z} } 
\nc{\tork}{ { T_k} }
\nc{\tilzed} {\tilde{Z} } 
\nc{\zed}{ {Z} }
\nc{\inclzed} {i_{\zed} }
\nc{\divis}[3]{D_{#1, #2}(#3)}
\nc{\tildediv}[3]{\tilde{D}_{#1, #2}(#3)}
\nc{\lamax}{ \bigwedge^{\rm max} }
\nc{\nminus}{{N^-}}
\nc{\omkd}{{ \omega_{n,d} }}
\nc{\secc}[2]{s_{#1}^{(#2)} }
\nc{\partlam}{\partial_\Lambda}
\nc{\partlamfun}[1]{\partial_{\Lambda_{#1}}}

\nc{\he}[1]{{\hat{u}_{#1} } }
\nc{\uroot}[1]{ {u_{#1}}  }
\newcommand{\labell}{\label}
\input amssym.def
\input amssym.tex
\newcommand{\renorm}{{ \setcounter{equation}{0} }}
\renewcommand{\theequation}{\thesection.\arabic{equation}}
\nc{\mnd}{{ M(n,d) } }
\nc{\weightl}{\Lambda^W}
\nc{\liet}{ {\bf t} } 
\nc{\lieg}{ { \bf g} } 
\nc{\ct}{ { \tilde{c} } }
\nc{\tilc }{ { \tilde{c} } }
\nc{\diag}{{ \rm diag} } 
\nc{\Res}{{\rm Res} } 
\nc{\res}{{\rm Res} } 
\nc{\lb}[1]{ {l_{#1} } }
\nc{\yy}[1]{Y_{#1}}
\nc{\liner}[1]{L_{#1} }
\nc{\lambdr}[1]{\Lambda_{#1} }
\nc{\linestd}{\liner{(\lb{1}, \dots, \lb{n-2} )} }
\nc{\lambstd}{\lambdr{(\lb{1}, \dots, \lb{n-2} )} }
\nc{\expsum}[1]{ { (e^{{#1} } - 1 ) } }
\nc{\itwopi}{ { 2 \pi i }}
\nc{\indset}{{(l_1, \dots, l_{n-1})}  }
\nc{\indsettwo}{{(l_1, \dots, \dots, l_{n-2})}  }
\nc{\bfl}{{ \bf \Lambda}}
\nc{\dgk}{{ D_{n,d}(g,k) }} 
\nc{\vgkl}{{ D_{n,d}(g,k, \Lambda) }} 
\nc{\vgklmult}{{ D_{n,d}(g,k, \bfl) }} 
\nc{\dgklmult}{{ V_{n,d}(g,k, \bfl) }} 
\nc{\vgk}{{ V_{n,d}(g,k) }} 
\nc{\dgkl}{{ V_{n,d}(g,k,\Lambda) }} 
\nc{\lineb}{L}
\nc{\resid}[1]{\res_{#1 = 0 } }
\nc{\residone}[1]{\res_{#1 = 1 } }
\nc{\kmod}{r}
\nc{\gmax}{{\g_{\rm max} }}
\nc{\constq}{ {\frac{1}{n}  }}
\nc{\sol}{ { S_{0 \mu} }}
\nc{\ch}{ { \rm ch} } 
\nc{\td}{ { \rm Td} } 
\nc{\abk}{\kappa}
\nc{\tg}{{\tilde{\gamma} } }
\nc{\inpr}[1]{ \langle #1 \rangle }
\nc{\ib}{\item{$\bullet$} }
\newcommand{\sg}{\Sigma^g} 
\title{The  Verlinde formula for  parabolic bundles}
\author{Lisa C. Jeffrey \\
Mathematics Department \\
University of Toronto\\
Toronto, Ontario M5S 3G3, Canada\\
alg-geom/0003150
\thanks{This material is based on work
supported by 
 grants from NSERC and the Alfred P. Sloan Foundation.
 \hspace*{\fill} MSC subject classification: 58F05} }
\maketitle
\begin{abstract}
Let $\sg$ be a compact Riemann surface of genus $g$, and
$G$ $= SU(n)$. 
We introduce the 
central element 
$c$ $ = \diag (e^{2 \pi i d/n}, \dots, e^{2 \pi i d/n} ) $ for 
$d$ coprime to  $n.$ 
In this paper we prove the Verlinde formula for 
the Riemann-Roch number of a line bundle over the 
moduli space $\modparonec$ of representations of the fundamental
group of a Riemann surface of genus $g$ with one boundary component,
for which the loop around the boundary  is constrained to 
lie in the conjugacy class of $c\exp (\Lambda)$ (for 
$\Lambda \in \lietpl) $,
and also for the moduli space $\modparc$ of representations
of the fundamental group of a Riemann surface of 
genus $g$ with $s+1$ boundary components for which
the loop around the 0-th boundary component is sent to the central 
element $c$ and 
the loop around the $j$-th boundary component
is constrained to lie in the conjugacy class
of $\exp (\Lambda^{(j)})$ for $\Lambda^{(j)} \in \lietpl$.
Our proof is valid for $\Lambda^{(j)} $ in  suitable
neighbourhoods of $0$.
\end{abstract}
\renorm
\section{Introduction}
Let $\sg$ be a compact Riemann surface of genus $g$, and
$G$ $= SU(n)$. 
We introduce the 
central element 
$c$ $ = \diag (e^{2 \pi i d/n}, \dots, e^{2 \pi i d/n} ) $ for 
$d$ coprime to  $n.$ 
In this paper we prove the Verlinde formula for 
the Riemann-Roch number of a line bundle over the 
moduli space $\modparonec$ of representations of the fundamental
group of a Riemann surface of genus $g$ with one boundary component,
for which the loop around the boundary  is constrained to 
lie in the conjugacy class of $c\exp (\Lambda)$ (for 
$\Lambda \in \lietpl) $. (Here, $\liet$ denotes the Lie algebra
of the maximal torus $T$ of $G$, and $\lietpl$ the fundamental 
Weyl chamber.)
We also prove the Verlinde formula
 for the moduli space $\modparc$ of representations
of the fundamental group of a Riemann surface of 
genus $g$ with $s+1$ boundary components for which
the loop around the 0-th boundary component is sent to the central 
element $c$ and 
the loop around the $j$-th boundary component
is constrained to lie in the conjugacy class
of $\exp (\Lambda^{(j)})$ for $\Lambda^{(j)} \in \lietpl$,
where we have introduced 
the notation 
$\bfl = (\Lambda^{(1)},\Lambda^{(2)}, \dots, 
\Lambda^{(b)}).$
The methods extend
the proof we gave in Section 11 of \cite{JK:mod} for 
the Verlinde formula for the  moduli 
space $\mnd$ of holomorphic vector bundles 
of coprime rank $n$ and degree $d$ and  fixed determinant, which
can alternatively be described as the space of representations
of the fundamental group of a Riemann surface of genus $g$ with
one boundary component into
$G$ which send the loop around the boundary to the central element
$c$.
Our proofs are  valid for $\Lambda$ in 
a suitable neighbourhood of $0$ and  $\bfl$  in
a suitable  neighbourhood
of ${\bf 0}$.

Our earlier work \cite{JK:mod} used the Riemann-Roch formula
to prove a formula for the dimension of the space of holomorphic
sections of powers of a certain line bundle over $\mnd$;
in the case of $\mnd$ the higher cohomology vanishes, so the 
Riemann-Roch formula gives the dimension of the zeroth cohomology
group. In this paper we exploit the fact that the  related
moduli spaces $\modparonec$ and more generally $\modparc$ defined above
(which
appear in algebraic geometry as 
moduli spaces of bundles with parabolic structure,
where the parameters
$\Lambda^{(s)}$ $\in \lietpl$ are equivalent 
to the specification of {\em weights};
see for example \cite{MS}) 
 fibre
over $\mnd$ provided the weights are sufficiently small. When these
spaces admit a prequantum line bundle, we may push forward along the 
fibre to obtain a formula for the Riemann-Roch number of the 
prequantum line bundle on the total space in terms of the evaluation
of appropriate cohomology classes
on the fundamental class of $\mnd$. 
Formulas for the intersection numbers in $\mnd$ were proved in 
\cite{JK:mod}; we apply these formulas together with
the fibration to recover the Verlinde
formula in this more general situation.

The layout of this paper is as follows. In Section 2 we review results from
\cite{JK:mod} on 
the cohomology ring of $\mnd$, while in Section 3 we summarize results
on symplectic fibrations. Section 4 contains our proofs of the
Verlinde formula: in Section 4.1 we first review the proof for 
$\mnd$ from \cite{JK:mod}, while Section 4.2 gives the proof for 
$\modparonec$ and Section 4.3 gives the proof for $\modparc$.

\noindent{\bf Acknowledgements:} This paper relies heavily 
on the author's earlier joint work with F. Kirwan, notably
on the paper \cite{JK:mod} where the Verlinde formula
for $\mnd$ is proved. We would like to acknowledge the
hospitality of 
Universit\'e Paris-Sud (Orsay), where part of the work was completed.

\renorm
\section{Review of results on the  cohomology  of moduli spaces}

\noindent{\bf Theorem} {\em (Atiyah-Bott 1982) } 
The cohomology ring of $\mnd$ over ${\bf Q}$ is  generated by elements
$\{ a_r, b_r^j, f_r \} $ (for $r = 2, \dots, n $ and 
$j = 1, \dots, 2g$), where 

$$ a_r \in H^{2r} (\mnd), $$
$$ b_r^j \in H^{2r-1} (\mnd), $$
$$ f_r \in H^{2r-2} (\mnd). $$

If $\univ$ is the universal vector bundle over $\mnd \times \sg$, 
the elements $a_r, b_r^j, f_r$ are obtained by decomposing 
$c_r(\univ)$ in terms of 
$$\oplus_{s = 0} ^2 H^{2r-s} (\mnd) \otimes H^s (\sg)$$
 using 
the K\"unneth formula. 
\begin{rem} The Verlinde formula (for $\mnd$ as well as for 
$\modparonec$ and $\modparc$) may be deduced from pairings 
involving only $f_2$ and the $a_j$.
\end{rem}

We review the results of \cite{JK:mod}, where formulas were proved
for intersection pairings in the cohomology of $\mnd$. 

\newcommand{\bracearg}[1]{ { [[ #1 ]] } }
\newcommand{\tildarg}[1]{ { [[ #1 ]]  } }

\begin{rem}
If $M$ is a symplectic manifold equipped with the Hamiltonian 
action of a Lie group $G$, and $0$ is a regular value of the moment
map $\mu: M \to {\bf g}^*$, then the Kirwan map is the map
$$\kappa: H^*_G(M) \to H^*_G(\mu^{-1}(0)) \cong H^*(M_{\rm red}).$$
In this paper we are using the notation $\kappa$  to refer to the
restriction to the image of $H^*_G({\rm pt}) \cong S({\bf g}^*)^G$
in the domain $H^*_G(M)$.
In the case where $M_{\rm red} = \mnd$,
this restricted map is simply the Chern-Weil map for the universal
bundle restricted to $\mnd \times \{ \rm pt \}.$
\end{rem}

\begin{theorem} \label{mainab}   Let $c=\diag\, (e^{2\pi i d/n},\ldots,
e^{2\pi i d/n})$ where $d \in\{1,\ldots,n-1\}$ is coprime to $n$,
and suppose that $\eta\in H^*_{G}(\{{\rm pt}\})$
is a polynomial $Q(\tau_2,\ldots,\tau_n)
$ in the 
equivariant cohomology classes $\tau_r \in H^{2r}_G({\rm pt}) $
$\cong S({\bf g}^*)^G$ 
 for $2\leq r\leq n$ (where $\tau_r$ is the $r$-th elementary
symmetric polynomial),
which generate the $G$-equivariant cohomology of a point, 
and map under the Kirwan map $\kappa$ to classes $a_r \in H^{2r}(\mnd)$. 
Let $f_2$ be the cohomology class of the symplectic form.
Then the pairing
$\kappa (\eta) \exp (f_2) [\mnd]$
is given by
$$\int_{\mnd} \kappa (\eta  )\exp (f_2) 
 =  \frac{(-1)^{n_+(g-1)}}{n!} \res_{Y_{1} = 0} \dots \res_{Y_{n-1} = 0} 
\Biggl ( \frac{\sum_{w \in W_{n-1}} e^{ 
\inpr{ \tildarg{w \tilc},X}  } \int_{T^{2g} } \eta
e^{ \omega} } {\cald^{2g-2} \prod_{1\leq j \leq n-1}
( \exp (Y_j)-1 ) } \Biggr ), $$           
where $n_+ = \frac{1}{2} n(n-1)$ is the number of positive
roots of $G=SU(n)$ and $X = (X_1. \dots, X_n)\in \liet \otimes \CC$ has 
coordinates $Y_1=X_1-X_2,\ldots,Y_{n-1}=X_{n-1}-X_n$ 
defined by
the simple roots, while $W_{n-1} \cong S_{n-1}$ is the Weyl
group of $SU(n-1)$ embedded in $SU(n)$ in the standard way
using the first $n-1$ coordinates. The quantity $\cald(X) = 
\prod_{\gamma> 0 } \gamma(X) $ is the product of the 
positive roots. The element
$\tilde{c}$  is the
unique element of $\liet$  which 
satisfies $\exp {\tilde{c}} = c$ and belongs to the
fundamental domain defined by the simple roots for the translation
action on $\liet$ of the integer lattice $\Lambda^I$ (in other words,
the fundamental alcove).

Also,  the notation $\bracearg{\gamma}$  
 means the unique element which is in the
fundamental domain defined by the simple roots for the
translation action  on $\liet$ of the integer lattice and for
 which $\bracearg{\gamma}$ is equal to $\gamma$ plus some element of the
integer lattice. 
\end{theorem}

Theorem \ref{mainab} is a special case of Theorem 
8.1 of \cite{JK:mod}: this special case suffices to prove the
Verlinde formulas.

\renorm
\section{Symplectic fibrations}

Let $\tilc \in \liet$ be an 
element of the closed fundamental alcove 
$D_+$ satisfying  $\exp \tilc = c$.
The interior of the fundamental alcove
will be denoted $\alco$.

Let $\sig{g}{n}$ denote an oriented 
 two-manifold of genus $g$ with $b$ oriented boundary
components $S_1, \dots, S_b$.

\begin{definition} \labell{upmodpardef}
Let ${\bt}  = (\orpar^{(1)}, \dots, \orpar^{(b)})
$ be a 
collection of $b$ values in $\alc$. The moduli
space of representations is defined by 
$$\modsppar{g}{b}(\bt)  = \upmodsppar{g}{b}(\bt)/\group,$$ where
\beq \labell{upmoddef} \upmodsppar{g}{b}(\bt) = 
 \{ \rho \in \Hom(\pi_1(\sig{g}{b}), \group): 
\rho([S_a] ) \in \conjclass (\exp \orpar^{(a)}) , a = 1, \dots, b\} \eeq
and $\group$  acts on $\upmodsppar{g}{b}(\bt) $ by conjugation. 
Here, $\conjclass (\exp \orpar^{[a]})$ denotes the conjugacy class 
of $\exp \orpar^{[a]}$ in $\group$. 
\end{definition}

The fundamental group of $\sig{g}{b} $ is  the free group
on $2g+b$ generators with one relation:
$$\pi_1(\sig{g}{n}) = < x_1, \dots, x_{2g}, y_1, \dots, y_b : 
\prod_{j = 1}^g [x_{j}, x_{j+g} ] = \prod_{r= 1}^b y_r > . $$
Thus we have 
\beq \labell{modpardef}
\modsppar{g}{n}(\bt) = 
\{ (h_1, \dots, h_{2g}, \beta_1, \dots, \beta_b) \in 
\group^{2g+b}: \prod_{j = 1}^g [h_{j}, h_{j+g} ] = 
 \prod_{r = 1}^b \beta_r, \beta_r \in \conjclass (\exp \orpar^{(r)})
\}/\group.  \eeq 
For later convenience we make a slight
modification of this definition:
We put 
\beq \labell{modpardefc}
\modparc = 
\{ (h_1, \dots, h_{2g}, \beta_1, \dots, \beta_b) \in 
\group^{2g+b}: \prod_{j = 1}^g [h_{j}, h_{j+g} ] = 
c \prod_{r = 1}^b \beta_r, \beta_r \in \conjclass (\exp \orpar^{(r)})
\}/\group,  \eeq 
where $c =  \diag (e^{2 \pi i d/n}, \dots, e^{2 \pi i d/n} ).  $
In particular we have
$$ \modparonec = 
\{ (h_1, \dots, h_{2g}, \beta) \in 
\group^{2g+1}: \prod_{j = 1}^g [h_{j}, h_{j+g} ] = c \beta, ~
\beta \in \conjclass (\exp \Lambda) \}/G. $$

Results on symplectic fibrations of moduli spaces were 
developed  
 in \cite{JW}, and used there for purposes distinct from the
objectives of the present paper.
 We have
\begin{theorem}   \labell{t:symformfibre}
There is a neighbourhood $U$ 
of $0$ in $\liet$
such that if $\orpar \in U$ then there is a fibration
\beq \label{orbitfibre} \pi: \modsppar{g}{1} (c,\orpar) \to \mnd \eeq
 with fibre 
$\calo_{\orpar} $ (the orbit of the adjoint action of $G$ on 
$\lieg$).
Further, the  symplectic form $\omega_\orpar$ on 
$\modsppar{g}{1}(\orpar)$ satisfies
\beq \label{symlam} \omega_{\orpar} = 
\pi^* \omkd + \tilde{\Omega}_{\orpar} \eeq
where 
$\omkd$ is the symplectic form on 
$\mnd$ and 
$\tilde{\Omega}_{\orpar} $ restricts on each fibre of $\pi$ to the 
standard Kirillov-Kostant symplectic form 
$\Omega_{\orpar} $ on 
the coadjoint orbit $\calo_{\orpar}$. 
\end{theorem}
\Proof This follows by general results regarding
symplectic fibrations associated to symplectic reduction at a regular value:
see \cite{J:ext}, Theorem 6.1 for a proof. Note that
$M (n,d)$ is obtained by reducing an appropriate extended
moduli space at a regular value of the moment map and that
the 2-form $\omega_\Lambda$ is nondegenerate in a neighbourhood of the
preimage of this regular value under the moment map:
see
Proposition 5.5 of \cite{J:ext}.  

\begin{rem} \labell{walldef} The region $U$ is characterized as the component
of the complement of the union 
of the walls $H_{v,n}$ in $\liet_+$ whose closure
contains $0$, where
$$H_{v,n} = \{ \Lambda \in \liet: v(\Lambda) = n \} $$
for one of the fundamental weights $v$ and $n \in \ZZ$. 
Clearly the $H_{v,n}$ are the hyperplanes where the function
\cite{qym} encoding the volume of the $\modparone$ is 
not  smooth  : in other words they bound  the region of 
regular values of the moment map, which is the region where
the symplectic fibration  of 
Theorem \ref{t:symformfibre} exists.
We thank A. Szenes for this observation, which
is explained in  \cite{Sz2}.
\end{rem}

\begin{prop} \labell{p:flagcoh}
 If $\orpar$ is a regular element of $\liets$, the coadjoint 
orbit
$\orb{\orpar}$ is diffeomorphic to the homogeneous space
$ \group/\torus$ 
so its cohomology is given by 
\beq 
H^*(\orb{\orpar}) \cong \frac{S(\liets)}{S(\liets)^W },  \eeq
in other words the quotient of  ring of polynomials on $\liet$ by 
the subring of symmetric polynomials. 
\end{prop}

We have the following proposition:

\begin{prop} \labell{p:flagmfd}
The space 
$\modparonec$ is  a splitting manifold for 
the universal bundle $\univ \mid_{\mnd \times {\rm \{pt\}} }$ over 
$\mnd \times {\rm \{pt\}} \subset \mnd \times \sig{g}{0}: $ 
in other words $$\pi^* \Biggl (\univ \mid_{\mnd \times {\rm \{pt\}} }\Biggr) = 
L_1 \oplus \dots \oplus L_n$$
 where $c_1(L_j) = e_j $ for a collection 
of classes $e_j $ in 
$H^2\Bigl (\modparonec
\Bigr )$ (for  $j = 1, \dots, n$). Here,  when 
$j = 1, \dots, n-1$, $e_j$ restricts on the 
fibres of $\pi$ to the generator
$\alpha_j$ $(j = 1, \dots, n-1)$ of $H^2 (G/T, \ZZ) 
\cong  H^1(T, \ZZ)$ 
corresponding to the $j$-th fundamental weight of $SU(n)$
(an element of ${\rm Hom} (T, U(1) ) $, which is 
isomorphic to $H^1 (T, \ZZ) $) and 
$e_n = - (e_1 + \dots + e_{n-1}). $ 
\end{prop} 
\Proof This follows from the algebro-geometric description
of the moduli space of parabolic bundles
(see for instance \cite{MS}): it is the moduli space parametrizing
 holomorphic
bundles over 
$\sig{g}{0}$ together with a flag in the fibre of each bundle over a basepoint
$(\{pt\}) \in \sig{g}{0}$. The flag structure enables us naturally to split
the universal bundle into a sum of holomorphic line bundles.

\begin{prop} \labell{p:spluniv} 
If  $\tau_r$ is the $r$-th elementary symmetric polynomial
(for $r = 2, \dots, n$) then $\tau_r(e_1, \dots, e_n) = \pi^* a_r$, 
where 
$a_r = c_r(\univ \mid_{\mnd \times {\rm  \{pt\} } } ). $
\end{prop}

\Proof See p. 284 (Section 21) of \cite{BT} for results on 
the properties of splitting manifolds and 
flag bundles. There, it is proved that for a complex vector bundle
$E$ over a complex manifold $M$ with splitting manifold 
$\flagm(E)$, we have 
\beq \labell{e:spluniv} H^* \Bigl (\flagm(E) \Bigr) = 
\frac{ H^*(M) [e_1, \dots, e_n] } { \prod_{i = 1}^n (1 + e_i) = c(E) }, \eeq
where the $e_j \in H^2(\flagm(E) ) $ restrict (for $j = 1, \dots, n$) 
on the fibre 
$U(n) /U(1)^n$ $\cong G/T$ (where $G  = SU(n)$ and 
$T$ is its maximal torus) 
to the images 
under the coboundary map (in the Leray-Serre spectral sequence) 
of the elements $H^1_{U(1)^n} (\{pt\}, \ZZ) $ 
$= \Hom(U(1)^{n}, U(1)) $ given by a basis for the weight lattice of 
$U(1)^n$.       

The following is a standard result (see for instance 
\cite{BGV}, Lemma 7.22):
\begin{prop} \labell{p:kks}
Let $\alpha_1, \dots, \alpha_n$ (subject to 
$\sum_{j = 1}^n \alpha_j = 0 $) be the basis for 
$H^2_{T} (\{pt\})  $ (the second equivariant cohomology group of a point
for the maximal torus $T$ of $SU(n)$) which
was introduced in Proposition 
\ref{p:flagmfd}.  Let $\Lambda = \sum_{i = 1}^{n-1} \Lambda_i \he{i}$
where $\he{i}$ is the $i$-th simple root: note that the simple roots
are the basis of 
$\liet$ dual to the fundamental weights, which were introduced in Proposition
\ref{p:flagmfd} to define the generators $\alpha_j$.
 Then the standard Kirillov-Kostant  symplectic 
form $\Omega_{\orpar } $ on 
$\orb{ \orpar } $ is given by 
$$ \Omega_{\orpar } = \sum_{j = 1}^n {\orpar}_j \alpha_j. $$
\end{prop} 

\renorm
\section{Application to the Verlinde formula}

\subsection{The Verlinde formula for $\mnd$}

\begin{definition}
The 
{\em  highest root} $ \gmax$ is given by 
$\gmax(X) = X_n - X_1$ or 
$\gmax(X)  = Y_1 + \dots + Y_{n-1} $. 
\end{definition}
\nc{\indek}{ { \Delta(r)} }

\begin{definition} \label{d10.1} The {\em Verlinde function}
$\vgk$ is given by 
$$\vgk = \sum_{\mu  \in \indek}
\frac{e^{-\itwopi \inpr{ \mu - \rho ,\tilc}  } } {(S_{0\mu}(k) )^{2g-2} } $$
where $\rho$ is half the sum of the positive roots and 
$$ S_{0\mu}(k) = \frac{1 }{\sqrt{n} \kmod^{(n-1)/2} }
\prod_{\g > 0 } 2 \sin \pi \inpr{\g,\mu} /\kmod . $$
\end{definition}
(See \cite{GW} (A.44) and \cite{qym} (3.16).)
Here, $\weightl$ is the weight lattice, identified
with points  in $\liet$.
We have introduced the quantity
\beq \labell{rdef}
\kmod = k + n; \eeq
we  have also introduced
$$ \indek = \{\mu \in 
\weightl_{\rm reg} \cap \lietpl : \inpr{\mu,\gmax} < r  \}. $$


\nc{\bracewtilc}[1]{ { \tildarg{w\tilde{c}}_{#1} } }
\nc{\bracewl}[1]{ { \tildarg{w (\tilc + v(\lambda + \rho)/r) }_{#1} } }
\nc{\bracewlmult}[1]{ {   
\tildarg{w (\tilc + \sum_{s = 1}^b v_s(\lambda^{(s)} + \rho)/r )  }_{#1} } }

\renewcommand{\ell}{{\call}}

The Verlinde formula is a formula for the dimension $\dgk$ 
of the space of holomorphic sections of 
powers of $\ell$ , where $\ell$ 
is a particular line bundle over $\mnd$: it has been proved
by Beauville and Laszlo \cite{BL}, Faltings \cite{F}, Kumar,
Narasimhan and Ramanathan \cite{KNR} and Tsuchiya, Ueno
and Yamada \cite{TUY}. In \cite{BiL} Bismut and Labourie
have given a proof of the Verlinde formula using
techniques from symplectic geometry.

In this section 
we review the results from 
Section 11 of \cite{JK:mod}, 
showing how the Verlinde formula follows from the formula
(Theorem \ref{mainab}) for intersection pairings in $\mnd$.

A line bundle $\ell$ over $\mnd$ may be defined for which
$c_1 (\ell) = n f_2$, since $n f_2 \in H^2(\mnd, \ZZ)$
(see \cite{DN}).  
 Whenever $k$ is a 
positive integer divisible by $n$, we then define
\beq \label{10.1} \dgk = \dim H^0(\mnd, \ell^{k/n} ). \eeq

Verlinde's conjecture 
 says that the Verlinde function 
specifies the dimension of the space of holomorphic sections
of $\ell^{k/n}$:
\begin{theorem}  \label{verl} {\bf (Verlinde's conjecture)} 
$$ \dgk = \vgk.$$  \end{theorem}

We review the method of Section 11 of 
\cite{JK:mod}, where we gave a proof of   Verlinde's conjecture for
$\mnd$:
 an outline  of the  method we use was given by 
Szenes \cite{Sz} (Section 4.2).

In fact  $H^i (\mnd, \ell^m ) = 0 $ for all $i > 0$ and
$m> 0$  (see Section 11 of \cite{JK:mod} for references and
an outline of the proof).
So
$\dgk$ is given for  $k>0 $ by the Riemann-Roch formula:

 \beq \label{10.2} \dgk = \int_{\mnd }\ch \ell^{k/n} \td \mnd. \eeq
We use the following results to convert (\ref{10.2}) into a form
to which we may apply our previous results.
\begin{lemma} \label{p10.1} 
For any complex manifold 
$M$ the Todd class of $M$ is given by 
$$\td(M) = e^{c_1(M)/2} \hat{A} (M) $$ 
where $c_1(M)$ is the first Chern class of the holomorphic tangent 
bundle of $M$, and $\hat{A}(M)$ is the $A$-hat genus of $M$.
\end{lemma}
\Proof See for example \cite{Gilkey}, pages 97-99.      

\begin{prop} \label{p10.2}
We have 
$$\hat{A} (\mnd) = \abk  \Bigl ( \prod_{\g > 0} 
\frac{ \g/2 }{\sinh \g/2 } \Bigr )^{2g-2}. $$
\end{prop}
\Proof This is proved by Newstead
in \cite{Newstead}.      
\begin{lemma} \label{p10.3} 
We have 
$$c_1(\mnd) = 2 n f_2. $$
\end{lemma}
\Proof This is proved in  \cite{DN}, 
Th\'eor\`eme F  .     

Of course the Chern character of $\ell^{k/n}$ is given by 
$\ch \ell^{k/n} = e^{k f_2}. $  
Thus we obtain
\begin{prop} \label{p10.4}
The quantity $\dgk$  is given by 
$$ \dgk = \int_{\mnd} e^{(k + n)f_2} \abk 
  \Bigl ( \prod_{\g > 0} 
\frac{ \g }{e^{\g/2}  - e^{-\g/2} } \Bigr )^{2g-2}. $$
\end{prop} 
\Proof This follows immediately from (\ref{10.2}), Lemmas 
\ref{p10.1} and \ref{p10.3} and Proposition \ref{p10.2}.

\begin{theorem} \label{t10.5} We have
$$\dgk = \frac{(-1)^{n_+(g-1)}}{n!}\sum_{w \in W_{n-1} } 
 \resid{Y_{1}}
\dots \resid{Y_{n-1} }  \Biggl ( e^{\kmod \inpr{\tildarg{w\tilc},X} } 
\int_{T^{2g}} e^{\kmod \omega} \times $$
\beq \prod_{\g > 0} 
\Bigl (\frac{ \g(X)}{e^{\g(X)/2} - e^{-\g(X)/2} }\Bigr )^{2g-2}  
\frac{1 } { \prod_{j = 1}^l (e^{ \kmod Y_j} - 1)
 \cald(X)^{2g-2} } \Biggr ).\eeq
\end{theorem}
\Proof This is a direct consequence of  Corollary 
\ref{p10.4} and Theorem \ref{mainab}.
Note that because the factor $e^{f_2}$ in the statement of 
Theorem \ref{mainab} has been replaced by $e^{\kmod f_2} $, it is necessary
to replace 
$  e^{\inpr{\tildarg{w\tilc},X} } $ by $e^{\kmod \inpr{\tildarg{w\tilc},X}
} $, 
and 
$e^{Y_j} - 1$ by $e^{\kmod Y_j} - 1$.

    We introduce $Z_j = \exp Y_j$. Since  for any $w \in W_{n-1}$ we have that
$$\tildarg{w\tilc} = \bracewtilc{1} \he{1} + 
\bracewtilc{2} \he{2} + \dots + \bracewtilc{n-1} 
\he{n-1}  $$ in terms of the standard basis vectors
$$\he{j}=(0,\ldots,0,1,-1,0,\ldots,0)$$ 
for the integer lattice $\intlat$ of $\liet$,
with  $n \bracewtilc{j} \in \ZZ $ 
for all $j$, and $0 \le \bracewtilc{j} < 1$ for 
all  $j$.
We obtain 
$$ e^{  \kmod \inpr{\bracewtilc,X} } = 
Z_1^{ \bracewtilc{1} \kmod } Z_2^{  \bracewtilc{2} \kmod } 
\dots Z_{n-1}^{ \bracewtilc{n-1} \kmod}. $$
(Recall that $k$ and $r$ are divisible by $n$ so 
$e^{  \kmod \inpr{\tilc,X} }  $ is a well defined single valued
function of $Z_1$, $\dots, Z_{n-1}$.)
Thus we can equate $\dgk$ with 
$$ 
 \frac{(-1)^{n_+(g-1) } }{n!}
\sum_{w \in W_{n-1} } \residone{Z_{1}} \dots \residone{Z_{n-1} }   
\Biggl ( \Bigl ( \prod_{j = 1}^{n-1} 
\frac{1}{Z_j}   \Bigr ) 
\int_{T^{2g}} e^{\kmod \omega} \times 
$$ 
$$
\frac{ Z_1^{ \bracewtilc{1} \kmod } Z_2^{  \bracewtilc{2} \kmod} \dots 
Z_{n-1}^{ \bracewtilc{n-1}\kmod} 
}{\prod_{\g > 0 } (\tg^{1/2}- \tg^{-1/2})^{2g-2} (Z_1^{  \kmod} - 1) 
\dots (Z_{n-1}^{ \kmod} - 1 ) } \Biggr )   $$ 
\beq \label{10.5} 
=  \frac{(-1)^{n-1 + n_+(g-1) } }{n!} \sum_{w \in W_{n-1} } 
\residone{Z_{1}} \dots \residone{Z_{n-1} }   
\Biggl ( \Bigl ( \prod_{j = 1}^{n-1} 
\frac{1}{Z_j}  \Bigr ) \times \eeq
$$ \int_{T^{2g}} e^{\kmod \omega} 
\frac{ Z_1^{ - \bracewtilc{1} \kmod}
 Z_2^{ - \bracewtilc{2}\kmod} \dots Z_{n-1}^{ -\bracewtilc{n-1} \kmod} 
}{\prod_{\g > 0 } (\tg^{1/2}- \tg^{-1/2})^{2g-2} (Z_1^{ - \kmod} - 1) 
\dots (Z_{n-1}^{ - \kmod} - 1 ) } \Biggr ). $$
Here, we have introduced $\tg  $ defined (for the root
$\g  = \uroot{r} + \uroot{r+ 1} + \dots + \uroot{s-1}$, where
the $\uroot{j} \in \liets$ are identified via the 
usual inner product with  the 
$\he{j} \in \liet$) by 
\beq \labell{tgdef} \tg(Z_1, \dots, Z_{n-1}) = Z_r \dots Z_{s-1}. \eeq
We also have 
\begin{lemma} \label{p10.6} 
$$\int_{T^{2}} e^{\omega}  =  n $$
and hence 
$$\int_{T^{2g}} e^{r \omega}  =  r^{(n-1)g} n^g. $$
\end{lemma}
\Proof This follows from the calculation
given  in Lemma 10.10 of \cite{JK:mod}.

The following may be proved by the same method as in Section 2 of 
\cite{JK:mod}
(see also \cite{Sz}):
\begin{prop} \label{p:vsz}  Suppose $\nu$ $\in \liet$ 
is of the form $\nu = \sum_{j = 1}^{n-1} \nu_j \he{j}$ with
$0 \le \nu_j <1$ for all $j$ (in other words 
$\nu$ is in the interior of the fundamental alcove). Define the 
meromorphic  function $f$ 
on the complexification
 $T^{\CC}$ of $T$ as follows: 
\beq \label{11.6} 
f(Z) = (-1)^{n-1} (-1)^{n_+(g-1)} \kmod^{(n-1)(g-1)} n^{g -1}
 \frac{Z_1^{-\nu_{1}\kmod} \dots Z_{n-1}^{-\nu_{n-1}\kmod} } 
{ \prod_{\g> 0 } (\tg^{1/2} - \tg^{-1/2} )^{2g-2} }.  \eeq
   Then we have  that
\beq\frac{1}{(n-1)!} \residone{Z_{1}} 
\dots \residone{Z_{n-1} } \sum_{w \in W_{n-1} } 
\prod_{j = 1}^{n-1} 
\Bigl ( \frac{\kmod }{Z_j} \Bigr )   \frac{\bracearg{ w f} (Z)} {\prod_{j = 1}^{n-1} 
(Z_j^{ - \kmod} - 1 ) } \eeq

$$= 
\sum_{\mu \in  \indek} 
f(\exp \itwopi \mu/\kmod). $$

Here, $W_{n-1} $ is the permutation group on 
$\{ 1, \dots, n-1 \} $ which is (isomorphic to) the Weyl group of
$SU(n-1)$, and $\bracearg{w f} $ is the function
\beq \label{11.6'} 
\bracearg{wf}(Z) = (-1)^{n-1} (-1)^{n_+(g-1)} \kmod^{(n-1)(g-1)} n^{g -1}
 \frac{Z_1^{-\bracearg{w\nu}_1\kmod } \dots 
Z_{n-1}^{-\bracearg{w\nu}_{n-1} \kmod}  } 
{ \prod_{\g> 0 } (\tg^{1/2} - \tg^{-1/2} )^{2g-2} }.  \eeq 
For a root $\gamma$, the quantity $\tg$ was defined by 
(\ref{tgdef}).

\end{prop}
\noindent{\em Remark:}
Notice that 
we have 
$$ \sum_{ \lambda \in \indek
  }  f(\exp \itwopi \lambda/r) = \frac{1}{n-1} 
\sum_{m_j = 1}^{r-1}  f\Bigl (e^{
\itwopi ( \sum_j m_j w_j )/r } \Bigr ) . $$ (Here,
the $w_j$ are the fundamental weights, which
are dual to the simple roots $\{ \he{j}\}$.) The set 
$\{ X \in  \liet: \; $ $ X = \sum_j \lambda_j \he{j}, 0 \le \lambda_j
< 1, \; j = 1, \dots, 
n-1 \} $
 is a 
fundamental domain for the action of the integer lattice $\Lambda^I$ 
on $\liet$, while the set 
$\{ X \in  \lietpl \subset \liet: \gamma_{\rm max} (X) < 1 \} $
 is a fundamental domain for the {\em affine Weyl group}
$W_{\rm aff}$ 
(the semidirect product of the Weyl group and the integer lattice),
and  $\intlat$
 has index $(n-1)! $ (rather 
than $n!$) in $W_{\rm aff} $
(in other words a fundamental domain 
for $\intlat$ contains $(n-1)!$ fundamental
domains for $W_{\rm aff}$). 

Applying  Proposition  \ref{p:vsz} 
we find (noting
that  $(-1)^{n-1} = c^\rho$ when $n$ and $d$ are coprime) that
\beq \label{10.10} 
\dgk =  (-1)^{n_+(g-1)} \kmod^{(n-1)(g-1)} n^{g -1} c^\rho 
\sum_{\l \in \indek}
\frac{ e^{ - \itwopi \inpr{\tilc, \l}  } }
 { \prod_{\g> 0 } (e^{\itwopi \inpr{\frac{\g}{2\kmod}, \l} } - 
e^{-\itwopi \inpr{\frac{\g}{2\kmod}, \l} } )^{2g-2}}. \eeq
This gives 
\beq \label{10.11}
\dgk = (-1)^{n_+(g-1)}  r^{(n-1)(g-1)} n^{g-1} 
\sum_{\l \in \indek}
\frac{ e^{ - \itwopi \inpr{\tilc, \l + \rho } }  }
{\prod_{\g > 0 } (2 i \sin \pi \inpr{\g,\l}/\kmod )^{2g-2} } \eeq
\beq \label{10.12}
 ~~~= \kmod^{(n-1)(g-1)} 
n^{g-1} 
\sum_{\l \in \indek}
\frac{ e^{ - \itwopi \inpr{\tilc, \l + \rho } }  }
{\prod_{\g > 0 }\Bigl (2  \sin 
\frac{\pi \inpr{\g,\l} }{\kmod }\Bigr )^{2g-2} }. \eeq
Comparing with Definition \ref{d10.1}, we see that $\dgk = \vgk$. 
This completes the proof of Theorem \ref{verl}.       


\subsection{The Verlinde formula for $\modparonec$}

There is a more general version of Verlinde's
conjecture which applies to the case of $\modparonec$.
 Let $\Lambda \in \liet$, and let $k$ be a positive integer divisible 
by $n$.  If $k \Lambda $ is in the weight lattice $\weightl$,
then the cohomology class of $k\omega $ 
is a class in integral cohomology, and hence is the first
Chern class of a line bundle over
$\modparonec$ (denoted $\call^{k/n}$).
Notice that $k f_2 $ is automatically in integral cohomology,
since $nf_2$ is in integral cohomology.
\begin{definition} 
$$\vgkl =  \sum_{j \ge 0 } (-1)^j {\rm dim} H^j (\modparonec, \call^{k/n})  $$
\end{definition}
Note that the argument sketched in Section 11 of
\cite{JK:mod} does not generalize to show that
$H^j(\modparonec, \call^{k/n}) = 0 $ when $j > 0 $, unlike
the situation for $\mnd$. (However, in \cite{Tel} Teleman 
has constructed an alternative argument showing the vanishing of these
higher cohomology groups.)

We introduce $\lambda = k \Lambda$.

\begin{definition}The Verlinde function is

$$\dgkl = 
\sum_{\mu \in \indek}
\frac{ e^{ - 2 \pi i \inpr{\mu - \rho, \tilc} } S_{\lambda \mu} (r) }
{ (S_{0 \mu} (r) )^{2g-1} } $$
where
\beq S_{\lambda \mu} (r) = 
\frac{ (\sqrt{-1})^{n(n-1)/2}  }{\sqrt{n} r^{(n-1)/2} } 
\sum_{w \in W}
 (-1)^w e^{ - 2 \pi i \inpr{w (\lambda+ \rho), \mu+ \rho}/r }. \eeq

\end{definition}

\begin{theorem}{\emph [Verlinde's conjecture for 
parabolic bundles] }\label{verlpar}
There exists a neighbourhood $U$ of $0$ in 
$\liet$ such that, for 
$\Lambda \in U$ for which
$k \Lambda \in \weightl$, we have
$$\vgkl  = \dgkl. $$
\end{theorem}

The proof of Theorem \ref{verlpar} proceeds
by a sequence of lemmas. By the Riemann-Roch formula,
 we have 
\beq \vgkl = \int_{\modparonec} {\rm ch} (\call^{k/n}) {\rm Td} (\modparonec)
\eeq
We now use the fibration  from Theorem \ref{t:symformfibre}
to integrate over the fibre in order to obtain 
an integral over the base $\mnd$: the latter is then
evaluated using
Theorem \ref{mainab}.
First we observe that because the Todd class is multiplicative,
it decomposes as the product of Todd classes corresponding to the
fibre and the base:
\begin{lemma}\labell{toddlem}
$${\rm Td} (\modparonec) = 
{\rm Td} (T_{\rm vert} \modparonec )
\pi^* {\rm Td} (\mnd). $$
Here we have introduced
$T_{\rm vert} \modparonec $, the
vertical tangent bundle of the fibre of
(\ref{orbitfibre}),
so
 $$\pi_* {\rm Td} (\modparonec ) = 
{\rm Td} (\mnd) \pi_* {\rm Td} (T_{\rm vert}\modparonec ). $$
\end{lemma}

Next we recall that since the fibre of (\ref{orbitfibre}) is
just a homogeneous space $G/T$, we can express its Todd
class in terms of the generators introduced in Proposition
\ref{p:flagmfd}.
\newcommand{\squige}{\cale}

\begin{lemma} \labell{l:toddcomp}
 $$ {\rm Td} (\calo_{\Lambda}) 
=  \left (\prod_{\gamma > 0 } \frac{\gamma(\squige) } {
1 - e^{- \gamma(\squige) } } \right )   .$$
Here, $\squige = (e_1, \dots, e_n)$ is regarded as a member 
of $H^2(\calo_{\Lambda }) \otimes \liet $ so one can naturally
pair it with the root $\gamma \in \liets$ to obtain an element 
$( \gamma, \squige)  = \gamma((e_1, \dots,
e_n))$ 
in $H^2 (\calo_{\Lambda },\RR).$
\end{lemma}
\Proof This is proved (for example) in  Section 14 of  \cite{Hirz}.
Here, the $e_i$ were introduced in Proposition \ref{p:flagmfd}.

To obtain the Riemann-Roch number we must compute
\beq \int_{\modparonec} {\rm ch} (\call^{k/n}) {\rm Td} (\modparonec). \eeq
Using the fact that
$$ {\rm ch} (\call^{k/n}) = \exp (k \omega_\Lambda) $$
and the decomposition of $\omega_\Lambda$ given by 
(\ref{symlam}), we find that we must compute
\beq \labell{integral}
\int_{\modparonec} \Biggl ( e^{k \pi^* \omega_{n,d} } \pi^* {\rm Td} (\mnd) 
\Biggr ) \Biggl ( 
e^{k \tilde{\Omega}_{\orpar } } 
{\rm Td} (T_{\rm vert} \modparonec ) \Biggr ) 
 \eeq
where $T_{\rm vert} \modparonec $ is the 
vertical tangent bundle of the fibres. 
\newcommand{\symnet}{{\calt(\cale)}}

\newcommand{\omnd}{\omega_{n,d}}
We shall integrate over the fibre appearing in (\ref{orbitfibre}):
on this fibre we see that
$$e^{k \tilde{\Omega}_{\orpar } }
{\rm Td} (T_{\rm vert} \modparonec )$$
becomes
\beq \labell{toddprod}\exp \left (\sum_{j=1}^n 
\lambda_j  e_j \right ) \prod_{\gamma >0}
\frac{ (\gamma, \squige) e^{ (\gamma, \squige)/2} }{ e^{(\gamma, \squige)/2}
- e^{-(\gamma, \squige)/2} } \eeq
Recall that  we have introduced
 $$\squige = (e_1, \dots, e_n).$$
To do the integral we reorganize it in terms of Weyl invariant
 and 
Weyl anti-invariant cohomology classes. 
We replace the expression in (\ref{toddprod}) by
\beq \labell{toddprodsym}
 \symnet    \prod_{\gamma >0 } (\gamma, \squige) , \eeq
where 
$$\symnet = \frac{1}{|W|} \sum_{v \in W} (-1)^v 
\frac{ e^{\inpr{v(\lambda  + \rho), \squige } }}{ \prod_{\gamma > 0 }
e^{(\gamma, \squige)/2} - e^{-(\gamma, \squige)/2}. } $$
Recall that $\rho$ is half the sum of the positive roots.

Note that the quantity $\symnet$ is
invariant under the transformation 
$\squige \mapsto w \squige$. Hence by Propositions \ref{p:flagcoh} and 
\ref{p:spluniv}, $\symnet$ $ = \pi^* \cals_\lambda(a_2, \dots, a_n)$ 
for an appropriate polynomial $\cals_\lambda(a_2, \dots, a_n) $
in the $a_2, \dots, a_n$.
Here, $\cals_\lambda$ is defined by 
\beq \labell{calseq} \cals_\lambda (\kappa(\tau_2), \dots, \kappa(\tau_n)) = 
\kappa \left(\frac{ \sum_{v \in W}  (-1)^v e^{(v(\lambda  + \rho), X) } }
{\prod_{\gamma > 0} e^{(\gamma, X)/2} - e^{- (\gamma, X)/2} } \right ). \eeq
The $a_j$ satisfy $a_j = \kappa (\tau_j)$ where
$\tau_j$ is the $j$-th elementary symmetric polynomial 
($\tau_j \in S(\liets)^W$) 
and $\kappa$ is the Kirwan map. The integral over the fibre 
reduces to
\beq \label{4.14a} \int_{\mnd}
\cals_\lambda(\kappa(\tau_2), \dots, 
\kappa (\tau_n) ) e^{k \omnd} {\rm Td}(\mnd) 
 \int_{\calo_{\Lambda }} \prod_{\gamma > 0 } (\gamma,\squige)
\eeq
where $\cals_\lambda$ is defined in (\ref{calseq}).
Now the integral $$\int_{\calo_{\Lambda }} 
\prod_{\gamma > 0} (\gamma, \squige) $$
gives the Euler characteristic of the orbit, which is simply
$|W|$ $= n!$ (\cite{Hirz}, Chap. 14). Thus the integral over 
the total space $\modparonec$ is reduced to an integral over
the base space $\mnd$: it becomes
\beq \labell{intlam}
\int_{\mnd} \cals_\lambda (\kappa (\tau_2), \dots \kappa(\tau_n))
e^{k \omnd}{\rm Td}(\mnd). \eeq
(Notice that the argument of 
$\kappa$ in (\ref{calseq}) is the quotient of two Weyl anti-invariant 
functions  of the variable $X$ and hence is Weyl invariant.)

A generalization of Theorem \ref{mainab} (proved by an extension of the
proof of Theorem \ref{mainab} given in \cite{JK:mod}) is

\newcommand{\posroots}{{\Delta_+}}

\begin{theorem} \label{mainabgen} 
In the notation from the statement of Theorem \ref{mainab},
let $$\alpha = \kappa\Biggl (
\sum_{v \in W} (-1)^v \frac{e^{\inpr{v \nu, X} }}{ 
\prod_{\gamma > 0} ( e^{\inpr{\gamma,X}} - e^{- \inpr{\gamma,X} } )} 
\Biggr ) $$
be a formal cohomology class in $H^*(\mnd)$,
where $\nu \in \liet^*$.  (Note that the argument of 
$\kappa$ is a Weyl invariant function of $X$.) Then 
$\kappa(\eta)\alpha \exp (f_2)\ [\mnd]$
is given by
 
$$\int_{\mnd} \kappa (\eta  )\alpha \exp (f_2)
 =  \frac{(-1)^{n_+(g-1)}}{n!} \res_{Y_{1} = 0} \dots \res_{Y_{n-1} = 0} 
$$
$$
\Biggl ( \frac{\sum_{w \in W_{n-1}} \sum_{v \in W} (-1)^v 
e^{\inpr{ \tildarg{w (\tilc+v \nu) } , X} } \int_{T^{2g} } \eta
e^{ \omega} } {\cald_n^{2g-2}
\prod_{\gamma > 0  } ( e^{\inpr{\gamma,X}/2} - e^{- \inpr{\gamma,X}/2 }  )
 \prod_{1\leq j \leq n-1}
( \exp (Y_j)-1 ) } \Biggr ), $$           
Recall that   the notation $\bracearg{\gamma}$  
 means the unique element which is in the
fundamental domain defined by the simple roots for the
translation action  on $\liet$ of the integer lattice and for
 which $\bracearg{\gamma}$ is equal to $\gamma$ plus some element of the
integer lattice. 
\end{theorem}

The expression in (\ref{intlam})
 is the type of integral computed by Theorem \ref{mainabgen}.
A straightforward modification of the proof of 
Proposition  \ref{p10.4} gives 
that the Riemann-Roch number is 
\beq 
\dgkl = \int_{\mnd} e^{(k + n)f_2} \abk 
  \Biggl ( \prod_{\g > 0} 
\frac { \g(X)^{2g-2} } { \left ( e^{\g(X)/2}  - e^{-\g(X)/2} \right )^{2g-1} }
\sum_{v \in W} (-1)^v e^{(v(\lambda + \rho ),X ) } \Biggr ). 
\eeq

Using Theorem \ref{mainabgen} 
this leads to (in a manner similar to the proof of
Theorem \ref{t10.5})
\beq \dgkl = 
 \frac{ (-1)^{n_+(g-1)} }{n!}\sum_{w \in W_{n-1} } 
 \resid{Y_{1}}
\dots \resid{Y_{n-1} }  
\int_{T^{2g}} e^{\kmod \omega} \times \eeq
$$ \prod_{\g > 0} 
\Bigl ( \frac{ (\g(X))^{2g-2} }{ (e^{\g(X)/2} - e^{-\g(X)/2})^{2g-1} 
}
\frac{\sum_{v \in W} (-1)^{v} e^{r\inpr{ \tildarg{ 
w (\tilc + v(\lambda + \rho)/r)},
X}  } }
 { \prod_{j = 1}^l (e^{ \kmod Y_j} - 1)
 \cald(X)^{2g-2} } .$$
This yields in turn (just as in (\ref{10.5}))
\beq \label{10.5pr} 
\dgkl=  \frac{(-1)^{n-1 + n_+(g-1) } }{n!} \sum_{w \in W_{n-1} } \sum_{v \in W}
(-1)^v\residone{Z_{1}} \dots \residone{Z_{n-1} }   
\Bigl ( \prod_{j = 1}^{n-1} 
\frac{1}{Z_j}  \Bigr ) \times \eeq
$$ \int_{T^{2g}} e^{\kmod \omega} 
\frac{ Z_1^{ - \kmod\bracewl{1} }
 Z_2^{ - \kmod\bracewl{2}} \dots Z_{n-1}^{ -\kmod\bracewl{n-1} } 
}{\prod_{\g > 0 } (\tg^{1/2}- \tg^{-1/2})^{2g-1} (Z_1^{ - \kmod} - 1) 
\dots (Z_{n-1}^{ - \kmod} - 1 ) }    $$  

Finally we may use Proposition \ref{p:vsz} to obtain an 
 expression equivalent to  (\ref{10.5pr}):
\begin{theorem} \label{t:dgkl} We have
$$
\dgkl = \sum_{\mu \in \indek}
\frac{ e^{ - 2 \pi i \inpr{\mu + \rho, \tilc} }  S_{\lambda \mu} (r) }
{ S_{0 \mu }^{2g-1} } $$
where
\beq S_{\lambda \mu} (r) = 
\frac{ (\sqrt{-1})^{n(n-1)/2}  }{\sqrt{n} r^{(n-1)/2} } 
\sum_{v \in W}
(-1)^v e^{ - 2 \pi i \inpr{v (\lambda+ \rho), \mu+ \rho}/r }. 
\eeq
\end{theorem}
Since the formula given in Theorem \ref{t:dgkl} is the formula
for $\vgkl$, we have proved Theorem \ref{verlpar}.
The proof of Theorem \ref{t:dgkl} is valid when $\Lambda \in U$.

\subsection{The Verlinde formula for $\modparc$}

\newcommand{\ga}[1]{{ \gamma^{(#1)} }}
More generally, there are spaces 
$\modparc$
which 
fibre  over $\mnd$ with fibre 
$\calo_{\Lambda^{(1)}  } \times \dots 
\times \calo_{\Lambda^{(b)}} $
for 
 $\Lambda^{(j)}$ sufficiently close to $0$.
We introduce the notation $\bfl = (\Lambda^{(1)}, \dots, \Lambda^{(b)})$.
We  restrict to those $\Lambda^{(j)}$ 
for which $k \Lambda^{(j)} \in \weightl$;
we denote $k \Lambda^{(j)}$ by $\lambda^{(j)}$.

\begin{prop} \label{p:fibmult}
There is a neighbourhood ${\bf U}$ of the element $(0, \dots, 0 )$ 
 $\in  \liet^b$ such that if
$\bfl \in {\bf U}$ then there is a 
fibration 
$$ \pi: \modparc \to \mnd$$ 
with fibre 
$\calo_{\Lambda^{(1)} } \times \dots 
\times \calo_{\Lambda^{(b)} } $. Further, the 
symplectic form on $\modparc$
 satisfies
\beq  \omega_{\bfl }
= \pi^* \omnd +  
\tilde{\Omega}_{\bfl } \eeq
where 
$\tilde{\Omega}_{\bfl} $ restricts on each fibre to the 
sum of the Kirillov-Kostant symplectic forms 
$$\sum_{j = 1}^b \Omega_{\Lambda^{(j)} } $$ on the 
product of the coadjoint orbits
$\calo_{\Lambda^{(j)}}. $
\end{prop}
\Proof The proof is as in Theorem \ref{t:symformfibre}. 
We use the fact that $\mnd$ is the reduced space at the
regular value ${\bf 0}$ of 
a symplectic manifold equipped with a Hamiltonian
action of $G^b$,  and that at orbits of the $G^b$ action
close to $0$, the corresponding reduced spaces are the
$\modparc$. This extends the argument of Proposition 5.5 of 
\cite{J:ext}, using the fact that the action of 
$G^b$ is free on the zero locus of the moment map: this
action is given in (5.9) and (5.10) of \cite{J:ext}, and it is 
straightforward to verify that the action is free when the moment
map takes the value $\{\bf 0\}.$ 
     
\begin{rem} \label{wallmult}
The region ${\bf U}$ is characterized as the 
component   whose closure contains $(0,\dots, 0)$
in the complement in ${\liet_+}^b$ 
of the union of the walls $$\Bigl \{ (\Lambda_1, \dots, \Lambda_b):
\sum_{(w_1, \dots, w_b) \in W^b} w_j \Lambda_j \in H_{v,n} \Bigr \} .$$
Here the walls $H_{v,n}$ were introduced in (\ref{walldef}).
This observation is  due to A. Szenes \cite{Sz2}.
\end{rem}

\begin{prop} \label{splitfieldtwo}
We have 
$$ \pi^* \Bigl [ \UU|_{\mnd \times {\{\rm pt\}} }\Bigr ]  = 
\oplus_{j = 1}^n \oplus_{s = 1}^b
L_j^{(s)}, $$
where $c_1(L_j^{(s)}) = e_j^{(s)}$ restricts on the $s$-th orbit 
$\calo_{\Lambda^{(s)}}$ in the fibre of
$\pi$ to the generator $\alpha_j^{(s)}$ of $H^2(\calo_{\Lambda^{(s)} })$
corresponding to the $j$-th fundamental weight of $SU(n)$.
\end{prop}

As in Lemma \ref{toddlem}, we have 
$${\rm Td} (\modparc) = 
{\rm Td} (T_{\rm vert} \modparc )
\pi^* {\rm Td} (\mnd), $$
where
${\rm Td} (T_{\rm vert} \modparc )$ is the vertical 
tangent bundle of $\modparc$.
The integral over the fibres decomposes as a product (as in 
(\ref{toddprod})) of factors of the form
\beq \labell{intgen}
\int_{\calo_{\Lambda^{(s)} } }\exp \left (\sum_{j=1}^b
\lambda^{(s)}_j  e_j^{(s)} \right ) \prod_{\gamma > 0 } 
\frac{ (\gamma, \squige^{(s)} ) e^{ (\gamma, \squige^{(s)})/2} }{ 
e^{(\gamma, \squige^{(s)} )/2}
- e^{-(\gamma, \squige^{(s)} )/2} } \eeq
where $\squige^{(s)} = (e_1^{(s)}, \dots, e_n^{(s)}). $
(Here we have introduced a basis $\{ e_j^{(s)} \} $ for 
$H^2(\calo_{\Lambda^{(s)}}  )$, following Propositions
\ref{p:flagcoh} and
\ref{p:flagmfd}.)
Each of the integrals (\ref{intgen}) is of the form
$$\cals_{\lambda^{(s)} }  (\kappa (\tau_2), \dots, \kappa(\tau_n))$$
where $\cals_\lambda$ was defined in (\ref{calseq}).
From the product of these integrals over the $\calo_{\Lambda^{(s)} }$,
we obtain
\beq \labell{4.23}
\dgklmult = \int_{\mnd}e^{r f_2} \times
\eeq $$ \times \kappa \Biggl ( 
 \prod_{\gamma>0} \frac {\gamma(X)^{2g-2} }{
\Bigl ( e^{\gamma(X)/2} - e^{-\gamma(X)/2} \Bigr )^{2g-2+b} }
\prod_{s = 1}^b \sum_{v_s \in W}   (-1)^{v_s}
 e^{ ( {v_s} (\lambda^{(s)}  + \rho), X)}
\Biggr ). $$

Using Theorem \ref{mainabgen}, (\ref{4.23}) leads to 
\beq \dgklmult = 
 \frac{ (-1)^{n_+(g-1)} }{n!}\sum_{w \in W_{n-1} } 
 \resid{Y_{1}}
\dots \resid{Y_{n-1} }  
\int_{T^{2g}} e^{\kmod \omega} \times \eeq
$$ \times \prod_{\g > 0} 
\Bigl ( \frac{ \g(X)^{2g-2} }{( e^{\g(X)/2} - 
e^{-\g(X)/2})^{2g-2+b} } \Bigr ) 
\prod_{s = 1}^b \sum_{v_s \in W}(-1)^{v_s} \frac{ 
\exp ( r \tildarg{ w \tilc+ \sum_s
(v_s(\lambda^{(s)} + \rho)/r},X)  }
 { \prod_{j = 1}^l (e^{ \kmod Y_j} - 1)
 \cald(X)^{2g-2} } .$$
This leads in turn (just as in (\ref{10.5}) and (\ref{10.5pr})) to

\beq \label{10.5pp} 
\dgklmult=  \frac{(-1)^{n-1 + n_+(g-1) } }{n!} \sum_{w \in W_{n-1} } 
 \sum_{v_1 \in W} \dots \sum_{v_b \in W}
(-1)^{v_1} \dots (-1)^{v_b}\times
\eeq 
$$\times \residone{Z_{1}} \dots \residone{Z_{n-1} }   
\Bigl ( \prod_{j = 1}^{n-1}   
 \frac{1}{Z_1} \times \dots \times 
\frac{1}{Z_{n-1} }  \Bigr ) 
\int_{T^{2g}} e^{\kmod \omega} $$
$$
\times \frac{ Z_1^{ - \kmod\bracewlmult{1} }
 Z_2^{ - \kmod\bracewlmult{2}} \dots Z_{n-1}^{ -\kmod\bracewlmult{n-1} } 
}{\prod_{\g > 0 } (\tg^{1/2}- \tg^{-1/2})^{2g-2+b} (Z_1^{ - \kmod} - 1) 
\dots (Z_{n-1}^{ - \kmod} - 1 ) }    $$  

Finally we may use Proposition \ref{p:vsz} to recover

\begin{theorem} \labell{rrmult}
The Riemann-Roch number of $\modparc$ is given by 
$$\int_{\modparc}
{\rm ch} (\call^{k/n}) {\rm Td} (\modparc) = 
 $$
$$  = 
\sum_\mu \frac{ e^{ - 2 \pi i \inpr{\mu - \rho, \tilc} }
\prod_{s = 1}^b 
S_{\mu \lambda^{(s)}} (r) }
{ (S_{0 \mu} (r) )^{2g-2+b} 
}. $$ 
\end{theorem}
The proof of Theorem \ref{rrmult} is valid when $\bfl
 \in   {\bf U}$.


\begin{thebibliography}{99}
\bibitem{AB} M.F. Atiyah, R. Bott, The Yang-Mills equations 
over  Riemann surfaces,  {\em Phil. Trans. Roy. Soc. Lond.} {\bf A308}
(1982) 523-615.
\bibitem{Beauv1} A. Beauville, Vector bundles on curves and 
generalized theta functions: recent results and open problems, 
{\em Current topics in complex algebraic
geometry (Berkeley, CA, 1992/93)}, 17--33, 
Math. Sci. Res. Inst. Publ., 28, 
Cambridge Univ. Press, Cambridge, 1995. 
\bibitem{BGV} N. Berline, E. Getzler, M. Vergne, 
{\em Heat Kernels  and Dirac Operators}, Springer-Verlag
(Grundlehren vol. 298), 1992.
\bibitem{BL} A. Beauville, Y. Laszlo, 
Conformal blocks and generalized theta-divisors, {\em Commun.
  Math. Phys.}  {\bf 164}  (1994) 385-519.
\bibitem{BiL} J.-M. Bismut, F. Labourie, Formule de Verlinde pour les groupes
simplement connexes et g\'eom\'etrie symplectique, {\em C. R. Acad.
 Sci. Paris S\'er. I Math.} 325 (1997), no. 9, 1009--1014; Symplectic 
geometry and the Verlinde formulas, Orsay preprint 98-66. 
\bibitem{BT} R. Bott, L. Tu, {\em Differential Forms in Algebraic 
Topology}  (Graduate Texts in Mathematics {\bf 82}). Springer-Verlag, 1982.
\bibitem{DN} J.-M. Drezet and M.S. Narasimhan, Groupe de Picard des
vari\'et\'es de modules de fibr\'es semi-stables sur les courbes
alg\'ebriques, {\em Invent. Math.} {\bf 97} (1989) 53-94. 
\bibitem{F} G. Faltings,
A proof of the Verlinde formula, 
{\em  J. Alg. Geom.}  {\bf 3} (1994) 347-374
\bibitem{GW} D. Gepner, E. Witten, String theory on group manifolds, 
{\em Nucl. Phys.} {\bf B 278} (1986) 493-549.
\bibitem{Gilkey} P.B. Gilkey, {\em Invariance Theory, The Heat 
Equation and the Atiyah-Singer Index Theorem}, Publish or Perish, 1984.
\bibitem{Hirz} F. Hirzebruch, \emph{ Topological Methods in 
Algebraic Geometry} (third edition),  Springer-Verlag, 1995.
\bibitem{J:ext} L.C. Jeffrey, Extended moduli spaces of 
flat connections on 
Riemann surfaces. {\em Math. Annalen} {\bf 298} (1994)    667-692.
\bibitem{JK:mod} L.C. Jeffrey and F.C. Kirwan, Intersection pairings
in moduli spaces of vector bundles of arbitrary rank on a 
Riemann surface.  {\em Annals of 
Math.} {\bf 148} (1998) 109-196.
\bibitem{JW} L.C. Jeffrey and J. Weitsman, Symplectic geometry of the
moduli space of flat connections on a Riemann surface: inductive
decompositions and vanishing theorems. {\em Canad. J. Math.},
to appear.
\bibitem{KNR} S. Kumar, M.S. Narasimhan,   A.  Ramanathan,
Infinite 
Grassmannian  and moduli spaces of $G$-bundles, 
{\em Math. Annalen} {\bf 300} (1994) 41-75.
\bibitem{MS} V.B. Mehta, C.S. Seshadri, Moduli of vector bundles on 
curves with parabolic structures. {\em Math. Ann.} {\bf 248} (1980),
205-239. 
\bibitem{Newstead} P. Newstead, Characteristic classes 
of stable bundles of rank 2 over an algebraic curve, 
{\em Trans. Amer. Math. Soc.} {\bf 169} (1972) 337-345. 
\bibitem{Sz} A. Szenes, The combinatorics of the Verlinde formula,
alg-geom/9402003,
in {Vector Bundles in Algebraic Geometry}, (Durham, 1993),
ed. N. Hitchin, W. Oxbury, Cambridge
University Press (LMS Lecture Series vol. 208), 
1995, 241-253; A. Szenes, private communication.
\bibitem{Sz2} A. Szenes, Iterated residues and 
multiple Bernoulli polynomials. 
{\em Internat. Math. Res. Notices} {\bf 18} (1998) 937--956. 
\bibitem{Tel} C. Teleman, The quantization conjecture
revisited, preprint math.AG/9808029.
\bibitem{TUY} 
A. Tsuchiya, K. Ueno, Y. Yamada, Conformal 
field theory on universal family of stable 
curves with gauge symmetry, {\em Advanced Studies in Pure Math.} 
vol. 19, Princeton Univ. Press, 1989, p. 459-566.
\bibitem{qym} E. Witten, On quantum gauge theories in two dimensions,
{\em Commun. Math. Phys.} {\bf 141} (1991) 153-209.
\bibitem{tdgr} E. Witten, { Two dimensional gauge theories
revisited}, preprint hep-th/9204083;
 {\em J. Geom. Phys.} {\bf 9} (1992) 303-368.
\end{thebibliography}
\end{document}